\documentclass[11pt]{amsart}

\makeatletter
\usepackage{amssymb}
\usepackage{latexsym}
\usepackage{amsbsy}
\usepackage{amsfonts}

\def\marginpar#1{\ignorespaces}

\newtheorem{theorem}[equation]{Theorem}
\newtheorem{proposition}[equation]{Proposition}

\newtheorem{definition}[equation]{Definition}
\newtheorem{assumption}[equation]{Assumption}

\theoremstyle{definition}
\newtheorem{remark}[equation]{Remark}

\newtheorem{example}[equation]{Example}

\numberwithin{equation}{section}

\def\AArm{\fam0 \rm}%
\newdimen\AAdi%
\newbox\AAbo%
\def\AAk#1#2{\setbox\AAbo=\hbox{#2}\AAdi=\wd\AAbo\kern#1\AAdi{}}%

\newcommand{\BBone}{{\ensuremath{{\AArm 1\AAk{-.8}{I}I}}}}

\def\eqref#1{(\ref{#1})}
\def\eqlabel#1{\def\@currentlabel{#1}}

\def\formula#1{\def\@tempa{#1}\let\@tempb\theequation\def\theequation{%
\hbox{#1}}\def\@currentlabel{(\theequation)}$$}
\def\endformula{\leqno\hbox{(\@tempa)}$$\@ignoretrue\let\theequation\@tempb}

\def\given{\hskip5\p@\relax\vrule\@width.4\p@\hskip5\p@\relax}

\newcommand{\open}[1]{%
\par\normalfont\topsep6\p@\@plus6\p@\trivlist\item[\hskip\labelsep\itshape#1%
\@addpunct{.}]\ignorespaces}

\DeclareRobustCommand{\close}[1]{%
  \ifmmode 
  \else \leavevmode\unskip\penalty9999 \hbox{}\nobreak\hfill
  \fi
  \quad\hbox{$#1$}}

\newlength{\toskip}

\def \R {{\mathbb R}}

\def \P {{\mathbb P}}
\def \E {{\mathbb E}}
\def \N {{\mathbb N}}

\def \L {{\mathbb L}}

\def \phi {\varphi}

\def \Var {\textrm{Var}}
\makeatother

\begin{document}
\date{\today}

\title[DEVIATION BOUNDS FOR ADDITIVE ...]{DEVIATION BOUNDS FOR ADDITIVE FUNCTIONALS OF MARKOV PROCESSES.}

 \author[P. Cattiaux]{\textbf{\quad {Patrick} Cattiaux $^{\spadesuit}$ \, \, }}

\address{{\bf {Patrick} CATTIAUX},\\ Ecole Polytechnique, CMAP, F- 91128 Palaiseau cedex,
CNRS 756\\ and Universit\'e Paris X Nanterre, \'equipe MODAL'X, UFR SEGMI\\ 200 avenue de la
R\'epublique, F- 92001 Nanterre, Cedex.} \email{cattiaux@cmapx.polytechnique.fr}

 \author[A. Guillin]{\textbf{\quad {Arnaud} Guillin $^{\diamondsuit}$}}
\address{{\bf {Arnaud} GUILLIN},\\ CEREMADE \, Universit\'e  Paris IX Dauphine, F- 75775 Paris cedex, CNRS
7534.} \email{guillin@ceremade.dauphine.fr}

\maketitle
 \begin{center}
 \textsc{$^{\spadesuit}$ Ecole Polytechnique \quad and \quad Universit\'e Paris X}
\medskip

\textsc{$^{\diamondsuit}$ Universit\'e Paris IX}
 \end{center}

\begin{abstract}
In this paper we derive non asymptotic deviation bounds for $$\P_\nu \left(\left|\frac 1t \,
\int_0^t V(X_s) ds \, - \, \int V d\mu \right| \, \geq \, R\right)$$ where $X$ is a $\mu$
stationary and ergodic Markov process and $V$ is some $\mu$ integrable function. These bounds are
obtained under various moments assumptions for $V$, and various regularity assumptions for $\mu$.
Regularity means here that $\mu$ may satisfy various functional inequalities (F-Sobolev,
generalized Poincar\'e etc...).
\end{abstract}

\bigskip

\bigskip

\bigskip

\section{\bf Introduction, framework and first results.}\label{Intro}

On some Polish space E, let us consider a conservative (continuous time) Markov process
$\left(X_t,(\P_x)_{x\in E}\right)$ and its associated semi-group $(P_t)_{t \geq 0}$ with
infinitesimal generator $L$ (and denote $D(L)$ its domain). Let $\mu$ be a probability measure on $E$ which is invariant and
ergodic w.r.t. $P_t$. The celebrated ergodic Theorem tells us that for any $V$ in $\L^1(\mu)$
$$A(t,R,V) \, := \, \P_\mu\left(\left|\frac 1t \, \int_0^t V(X_s) ds \, - \, \int V d\mu \right| \,
\geq \, R\right) \, \to \, 0$$ as $t$ goes to $+\infty$ for all $R>0$. Level 1 large deviations
theory furnishes asymptotic bounds for $\frac 1t \, \log(A(t,R,V))$ (see e.g. \cite{DS}). If $V$ is
bounded, one may replace the initial measure $\mu$ by a $\mu$ absolutely continuous probability measure $\nu$.

It is however of major importance in practice to exhibit non asymptotic upper bound but also to ensure practical conditions to verify them; see for example a priori bounds for large and moderate deviations in averaging principle, concentration for particular approximations of granular media equation,... . It will be the purpose of the present note. In \cite{Wu0}, Liming Wu derived such bounds. The main result by Wu
reads as follows: if $V$ is bounded, then for all $t>0$ and all $R>0$
\begin{equation}\label{eqwu1}
\P_\nu\left(\frac 1t \, \int_0^t V(X_s) ds \, - \, \int V d\mu  \, \geq \, R\right) \, \leq \,
\parallel \frac{d\nu}{d\mu}\parallel_{\L^2(\mu)} \, \exp \left\{ - \, t \, I_V\left(R+\int V d\mu\right) \right\}
\end{equation}
where $I_V(a) \, = \, \sup_{\lambda \geq 0} \{\lambda \, a \, - \, \Lambda(\lambda V)\}$ and
$$\Lambda(V) := \sup \, \left\{ \int \, Vf^2 \, d\mu \, + \, <Lf,f>_\mu \, ; \, f\in D(L) \,
\textrm{ and } \, \int f^2 d\mu =1 \right\} \, .$$ Of course a similar result holds for
$\P_\nu\left(\frac 1t \, \int_0^t V(X_s) ds \, - \, \int V d\mu  \, \leq \, -R\right)$.

The key is that $$\frac 1t \, \log \parallel P_t^V\parallel_{\L^2(\mu)} \, \leq \, \Lambda(V)$$ for
all $t>0$, where $P_t^V$ denotes the Feynman-Kac semi-group built from $P_t$. This result is a
consequence of Lumer-Philips Theorem. It is worthwhile noticing that, when $P_t$ is $\mu$
symmetric, the above bound is asymptotically sharp, according to the spectral radius theorem, but
\eqref{eqwu1} is also asymptotically sharp according to large deviations theory (see \cite{DS}
Theorem 5.3.10).
\smallskip

The main difficulty is then to be able to give a precise (and if possible optimal) control of the quantity $I_V(a)$ and by way of $\Lambda(\lambda V)$. Our approach mainly relies on the use of functional inequalities to get upper bound on $\Lambda(\lambda V)$. Let us illustrate this approach via the use of a Poincar\'e inequality (or spectral gap inequality).
\smallskip

Take first $V$ bounded. Of course, by homogeneity, we may only consider the $V$'s satisfying $\int V d\mu=0$ and $\sup |V|=
1$, for which the only interesting $R$'s are between 0 and 1. Indeed notice that the bound
\eqref{eqwu1} is fortunately 0 if $R>1$ in this case since $\Lambda(\lambda V) \leq \lambda$ so
that $I_V(R)=+\infty$ if $R>1$. The next result furnishes an explicit bound as soon as a $\mu$
satisfies a Poincar\'e inequality

\begin{proposition}\label{propPoinc}
Assume that $\mu$ satisfies the Poincar\'e inequality $$\Var_{\mu}(f) \leq - C_P \, <Lf,f>_\mu \,
.$$ Then for all $V$ such that $\sup |V|= 1$, all $0<R\leq 1$ and all $t>0$
\begin{equation}\label{eqlez1}
\P_\nu\left(\frac 1t \, \int_0^t V(X_s) ds \, - \, \int V d\mu  \,
\geq \, R\right) \, \leq \,
\parallel \frac{d\nu}{d\mu}\parallel_{\L^2(\mu)} \, \exp \left( - \, \frac{t R^2}{8 C_P \Var_\mu(V)}\right) \, .
\end{equation}
\end{proposition}

\begin{proof}
We may assume that $\int V d\mu =0$. If $\int f^2 d\mu =1$ we may write $$f=\frac{1+\varepsilon
g}{\sqrt{1+\varepsilon^2}}$$ for some $\varepsilon \geq 0$ and some $g$ satisfying $\int g d\mu =0$
and $\int g^2 d\mu =1$, and conversely. Thus applying Poincar\'e with $G_P=1/C_P$
\begin{eqnarray*}
\Lambda(\lambda V) & \leq & \sup \, \left\{ \int \, \lambda Vf^2 \, d\mu \, - \, G_P \,
\Var_{\mu}(f) \, ; \, f\in D(L) \, \textrm{ and } \, \int f^2 d\mu =1 \right\} \, \\ & \leq &
\sup_{\varepsilon \geq 0} \, \left(\frac{\varepsilon}{1+\varepsilon^2} \, \sup \, \left\{2\lambda
\, \int Vg d\mu + \varepsilon \, \int (\lambda V - G_P) g^2 d\mu \, \right\}\right)
\end{eqnarray*}
where the second supremum is taken over the set $\left\{ g\in D(L) \, , \, \int g^2 d\mu =1 \, , \,
\int g d\mu =0 \, \right\}$.

It follows according to Cauchy-Schwarz and our hypotheses $$\Lambda(\lambda V) \leq
\sup_{\varepsilon \geq 0} \, \frac{\varepsilon}{1+\varepsilon^2} \, \left(2\lambda  \, \Var_\mu(V)
\, + \, \varepsilon \, (\lambda - G_P)\right) \, ,$$ so that, bounding roughly
$1/(1+\varepsilon^2)$ by $1$, we finally obtain that for $\lambda < G_P$, $$\Lambda(\lambda V) \leq
\frac{\lambda^2 \, \Var_\mu(V)}{G_P - \lambda} \, .$$ Taking the supremum on $\{ \lambda \in
[0,G_P/2] \}$ we finally deduce that
$$I_V(R) \geq \frac{G_P R^2}{8 \Var_\mu(V)} \, .$$
\end{proof}

We did not try to obtain the sharpest bounds in the previous Proposition because a similar (and
a little more precise) result was obtained by Lezaud (\cite{Lez}) using Kato's perturbation theory.
Actually the best uniform result contained in \cite{Lez} is exactly ours, but Lezaud obtains very
interesting non uniform results. Our proof above is much shorter.
\medskip

An interesting feature is that Proposition \ref{propPoinc} admits a (partial) converse. Indeed
\begin{proposition}\label{propPoinc2}
Assume that $\mu$ is diffuse (i.e. for any $A$ and any $0 \leq \alpha \leq \mu(A)$ there exists
$B\subseteq A$ such that $\mu(B)=\alpha$).

Assume that there exist $C$ and $\lambda_0$ such that for all $V$ such that $\int V d\mu=0$ and
$\sup |V|= 1$ and all $0 \leq \lambda \leq \lambda_0$, $\Lambda(\lambda V)\leq C \lambda^2$. Then
$\mu$ satisfies a Poincar\'e inequality.
\end{proposition}
\begin{proof}
Using the same notation as before, we deduce from the hypotheses that for all $\varepsilon \geq 0$
and all $g$ such that $\int g d\mu=0$ and $\int g^2 d\mu =1$, for $0 \leq \lambda \leq \lambda_0$,
\begin{equation*}
 0 \leq  C(1+\varepsilon^2) \lambda^2 -  \lambda \, \left( \varepsilon^2 \int Vg^2 d\mu + 2 \varepsilon \int Vg
 d\mu\right)\,  - \,  \varepsilon^2 \, <Lg,g>_\mu \, .
\end{equation*}
Notice that the above quantity reaches its minimum for $$\lambda = \frac{2\varepsilon \int Vg d\mu
+ \varepsilon^2 \int Vg^2 d\mu}{2C(1+\varepsilon^2)}$$ that goes to $0$ when $\varepsilon$ goes to
$0$ and to $\int  Vg^2 d\mu / 2C \, \leq 1/2C$ when $\varepsilon$ goes to $+\infty$. Thus, taking a
larger $C$ if necessary, we may assume that $1/2C \leq \lambda_0$, and changing $V$ into $-V$ if
necessary, we may assume that $\int V g^2 d\mu \geq 0$.

For $\varepsilon$ small enough the minimum is reached at some $\lambda \leq \lambda_0$ and has to
be nonnegative. It follows
$$\left(2\int Vg d\mu + \varepsilon \int Vg^2 d\mu\right)^2 \leq 4C (1+\varepsilon^2) <-Lg,g>_\mu
\, ,$$ so that letting $\varepsilon$ go to $0$ we obtain
$$\left(\int Vg d\mu\right)^2 \, \leq \, - \, C <Lg,g>_\mu \, .$$ We may then choose
$V=\textrm{sign}(g) - \int \textrm{sign}(g) d\mu$ in order to obtain
\begin{equation}\label{eqmoy}
\left(\int |g| d\mu \right)^2 \, \leq \, - \, C <Lg,g>_\mu \, .
\end{equation}

For $\varepsilon$ going to $+\infty$ and provided $\int V g^2 d\mu \geq 0$ we also obtain
\begin{equation}\label{eqvgcarre}
\left(\int V g^2 d\mu \right)^2 \, \leq \, - \, 4C <Lg,g>_\mu \, .
\end{equation}
We shall now build an appropriate $V$.

Let $A=\{|g|\leq 1/2\}$. First, if $\mu(A) \, \leq \, 1/2$, $\int |g| d\mu \, \geq \, 1/4$, so that
\eqref{eqmoy} implies $$\int g^2 d\mu = 1 \leq - \, 16 C <Lg,g>_\mu \, .$$ If $\mu(A) \, \geq \,
1/2$, denote by $B=\{|g|\geq 3/4\}$. We have $$1 = \int g^2 d\mu \leq \int_{B^c} g^2 d\mu + \int_B
g^2 d\mu \leq \frac{9}{16} + \int_B g^2 d\mu$$ so that $ \int_B g^2 d\mu \geq 7/16$. Choose
$A'\subseteq A$ such that $\mu(A')=\mu(B)$, and $V=\BBone_B - \BBone_{A'}$. Then $$\int V g^2 d\mu
\geq \frac{7}{16} - \frac 14 \, \mu(B) \, \geq \, \frac{3}{16} \geq 0 \, ,$$ and $$\int g^2 d\mu =
1 \, \leq \, \left(\frac{16}{3}\right)^2 \, 4C <-Lg,g>_\mu \, .$$ Hence there exists some constant
$K$ such that $\int g^2 d\mu \leq - \, K <Lg,g>_\mu$ for all $g$ with mean 0 and variance 1, that
is Poincar\'e holds.
\end{proof}

\begin{remark}
It is easy to see that \eqref{eqmoy} (that holds without the assumption of diffusivity) implies the
following
$$\Var_\mu^2(g) \, \leq \, C \, <-Lg,g>_\mu \, \,
\parallel g \parallel_\infty^2 \, ,$$
which is some weak Poincar\'e inequality.  This inequality implies some concentration property for
$\mu$ (see e.g. \cite{r-w} or \cite{BCR2}) but is quite far from the usual Poincar\'e inequality.
More precisely the aforementioned weak Poincar\'e inequality on $\R$ implies that $\mu$
concentrates like $\alpha(ds) = c/(1+|s|^3) \, ds$ and is actually satisfied by $\alpha$.
\end{remark}
\smallskip

It is not difficult to see that $\Lambda(\lambda V)\leq C \lambda^p$ for some $p>2$ (and small
$\lambda$'s) cannot happen (using the same method). This is natural since for very small $R$ we
cannot expect a better behaviour as a Gaussian one, due to the Central Limit Theorem (see e.g.
\cite{Lez} Theorem 3.1).
\medskip

We can now state the problems we shall study in the sequel :
\begin{itemize}
\item What happens if Poincar\'e is reinforced, replacing it by stronger functional inequalities~?
The answer to this question is partly given in \cite{Wu0} for the log-Sobolev inequality, namely we
may consider in this case unbounded $V$ having some exponential moments. \item What can be said for
bounded $V$'s when Poincar\'e's inequality does not hold ? \item For unbounded $V$'s, how can we
obtain (may be rough) deviation bounds in full generality ? \item What happens if the initial
measure is no more absolutely continuous, or when its density is less integrable ?
\end{itemize}
\medskip

\section{\bf Exponential bounds for unbounded $V$'s and strong functional inequalities}\label{sob}

Let us start here with an almost immediate extension of Wu's
result, tackling the first question.

\begin{theorem}\label{thmfsob}
Let $F$ be defined on $\R^+$. We assume that $F$ is continuous, increasing, concave, goes to
$+\infty$ at $\infty$ and satisfies $F(1)=0$. It follows that $F$ admits an inverse function which
is defined on $]F(0),+\infty[$. In addition we assume that $F$ satisfies
\begin{equation}\label{eqsub}
F(xy) \, \leq  \, F(x) \, + \, F(y) \, ,
\end{equation}
for all positive $x$ and $y$. If $\mu$ satisfies the following $F$-Sobolev inequality
$$\int f^2 \, F(f^2) \, d\mu \, \leq \, -<Lf,f>_\mu \, ,$$ for all $f \in D(L)$ such that $\int f^2
d\mu =1$, then for all $R>0$
\begin{equation}\label{eqfsob}
\P_\nu\left(\frac 1t \, \int_0^t V(X_s) ds \, - \, \int V d\mu  \, \geq \, R\right) \, \leq \,
\parallel \frac{d\nu}{d\mu}\parallel_{\L^2(\mu)} \, \exp \left( - t \, H^*(R+\int V d\mu)\right) \, ,
\end{equation}
where $$H^*(a) := \sup_{0\leq \lambda < \lambda_V} \, \left\{\lambda a \, - \, F\left(\int
F^{-1}(\lambda V) \, d\mu \right)\right\} \, ,$$ where $\lambda_V$ is such that $\lambda V > F(0)$
for all $0 \leq \lambda < \lambda_V$.

We also have $$\P_\nu\left(\frac 1t \, \int_0^t V(X_s) ds \, - \, \int V d\mu  \, \geq \, R\right)
\, \leq \,
\parallel \frac{d\nu}{d\mu}\parallel_{\L^2(\mu)} \, \exp \left( - t \, H_c^*(R)\right) \,
,$$ with $$H_c^*(a) := \sup_{0\leq \lambda < \lambda'_V} \, \left\{\lambda a \, - \, F\left(\int
F^{-1}(\lambda (V - \int V d\mu)) \, d\mu \right)\right\} \, ,$$ where $\lambda'_V$ is such that
$\lambda (V - \int V d\mu) > F(0)$ for all $0 \leq \lambda < \lambda'_V$. This latter bound is
better than the previous one when $F(0)=-\infty$.
\end{theorem}
\begin{proof}
Assume first that $V$ is bounded. Applying the $F$-Sobolev inequality, we get
\begin{eqnarray*}
\Lambda(\lambda V) & \leq & \sup \, \left\{ \int \, \lambda V \, f^2 \, d\mu \, - \, \int f^2
F(f^2) d\mu \, ; \, f\in D(L) \, \textrm{ and } \, \int f^2 d\mu =1 \right\} \, ,
\end{eqnarray*}
so that for all $f$ as above
\begin{eqnarray*}
F\left(\int F^{-1}(V) d\mu\right) & = & F\left(\int F^{-1}(V) d\mu\right)\\
& = & F\left(\int \frac{F^{-1}(V)}{f^2} \, f^2 d\mu\right)\\ & \geq & \int \,
F\left(\frac{F^{-1}(V)}{f^2}\right) \, f^2 d\mu \\ & \geq & \int \left(V \, - \, F(f^2)\right) f^2
d\mu \,
\end{eqnarray*}
where we have successively used the facts that $F$ is non-decreasing, concave and \eqref{eqsub}. It
follows that
$$\Lambda(\lambda V) \, \leq \, F\left(\int F^{-1}(\lambda V) d\mu\right) \, .$$ If $V$ is not
bounded just approximate it by $(V\wedge n)\vee -n$.

Finally we may replace $V$ by $V - \int V d\mu$ and obtain the last statement. The property $F(xy)
\leq F(x)+F(y)$ immediately shows that this bound is better than the previous one except that the
authorised set of $\lambda$'s differ in general.
\end{proof}

\begin{remark}
The bound obtained in Theorem \ref{thmfsob} is interesting since, assuming some regularity for $F$,
$\Lambda(\lambda (V - \int V d\mu))$ behaves like $\lambda^2$ for small $\lambda$ provided it is
finite for some $\lambda_0 >0$. Hence $H_c^*$ is strictly positive and actually behaves like $C
a^2$ for small $a$ while it behaves like $C a$ for large $a$.

Note that if $F(0)>-\infty$ the Theorem only applies to the bounded from below $V$'s.
\end{remark}

In the examples below $-<Lf,f>_\mu = 1/2 \, \int |\nabla f|^2 d\mu$, corresponding to diffusion process with constant diffusion term.
\begin{example}\label{exF}
1) The function $F(x)=C \, \log(x)$ satisfies all the previous assumptions with $F(0)=-\infty$. The
corresponding result is then Corollary 4 in \cite{Wu0}. Gaussian measures satisfy such log-Sobolev
inequalities. In this case $F^{-1}(y)=\exp (y/C)$, so that the above result holds as soon as $V$
has some exponential moment.

Note that some converse holds in this case. Indeed if $F= C \log$, and if $$\Lambda(V) \, \leq \,
\frac 1C \, \log\left(\int e^{CV} d\mu\right) \, ,$$ for all $V$, then for all $f$ we may choose
$V=\frac 1C \, \log f^2$ and deduce the log-Sobolev inequality.
\medskip

 2) The functions $F_\alpha(x) = \log^{\alpha}(1+x) \, - \, \log^{\alpha}(2)$ also satisfy all
the assumptions as soon as $0<\alpha\leq 1$ (see the proof of Theorem 38 in \cite{BCR1} for
instance). The measure $\mu_\beta (dx) = \frac{\exp ( - |x|^\beta)}{Z_\beta} \, dx$ satisfies a
 $c_\alpha F_\alpha$-Sobolev inequality for $\alpha = 2(1-1/\beta)$ and some well
chosen constant $c_\alpha$ (see \cite{BCR1} section 7). Here $F_\alpha^{-1}(y)$ behaves like $\exp
(y^{1/\alpha})$ at infinity. Here again some converse holds, but details are a little bit tedious.
\medskip

3) Conditions for some $F$-Sobolev inequalities are discussed in details in \cite{BCR1} and
\cite{BCR3}. In particular explicit (and tractable) criteria for absolutely continuous measures on
the line are given in \cite{BCR1} Theorem 27, while sufficient conditions are discussed in
\cite{BCR3} section 8 for a general Riemannian manifold. In these papers the corresponding
$F$-Sobolev inequalities satisfy the tensorization property, hence due to the CLT, $F$ cannot grow
faster than a logarithm. That is the described field of measures is between Poincar\'e and Gross
(log-Sobolev) inequalities. The condition $F(xy)\leq F(x)+F(y)$ certainly obliges us to such a
restriction.
\end{example}
\smallskip

\begin{example}\label{exgauss}
It is interesting to see how the previous result applies on simple examples. Consider the standard
Ornstein-Uhlenbeck process on $\R$, $dX_t=dB_t \, - \, \frac 12 \, X_t dt$ with its symmetric
probability measure $\gamma$ the standard normal law. $\gamma$ satisfies a log-Sobolev inequality
with constant $C=4$. Easy calculations yield
\begin{equation*}
\P_\nu\left(\frac 1t \, \int_0^t X_s^2 ds \, - \, 1  \, \geq \, R\right) \, \leq \,
\parallel \frac{d\nu}{d\mu}\parallel_{\L^2(\mu)} \, \exp \left( - \frac t8 \, (R - \log(1+R)) \right) \,
.
\end{equation*}
This result is not asymptotically sharp, since according to a result by Bryc and Dembo
(\cite{Bryc}) the large deviations rate function is $R^2/8(R+1)$ which is greater than our $(R -
\log(1+R))/8$. In particular for small $R$, we are loosing a factor 2. Note that Lezaud obtains in
\cite{Lez} Example 4.2 the correct bound, but that this case is a little bit miraculous since the
spectral gap of the Feynman-Kac operator can be explicitly calculated. We shall discuss other
explicit examples later on.
\end{example}
\smallskip

\begin{remark}
In the examples above, we have assumed that the diffusion coefficient is constant. And one knows that $F$-Sobolev inequalities are usually verified with the energy given by $\int |\nabla f|^2$ which could be seen as a limitation on the diffusion process we may consider for our deviation inequalities. However we may easily replace this assumption by some strict ellipticity, namely suppose that there exists $\delta>0$ such that for all $x,y$, $\langle \sigma(x)\sigma(x)^*y,y\rangle \ge \delta |y|^2$, then$$-\langle Lf,f\rangle={1\over 2}\int|\sigma(x)\nabla f(x)|^2d\mu(x)\ge{\delta\over 2}\int|\nabla f|^2d\mu.$$ 
It enables us to consider deviation inequalities for strictly elliptic diffusion of the form
$$dX_t=b(X_t)dt+\sigma(X_t)dW_t$$
using ``standard'' functional inequalities.
\end{remark}

\begin{remark}
One strongly suspects that the integrability condition ($\int F^{-1}(\lambda V) d\mu < +\infty$ for
some $\lambda >0$) is also necessary for an exponential bound to hold. We do not know whether this
is true in full generality or not, but one can easily build some examples.

Still in the Gaussian case of example \ref{exgauss}, consider $V(x)=x^4$, and define $V_N=V\wedge
N$. Choosing $f(x)=c \, e^{x^2/4}/(1+x^2)$ for some normalising constant $c$, we immediately see
that for $N$ large and $\lambda \geq 1/N^{\frac 18}$, $$\Lambda(\lambda V_N) \geq D \lambda \,
N^{\frac 14} \geq D N^{\frac 18}$$ for some nonnegative constant $D$. It follows that $H^*_N(a) \,
\leq \, a/N^{\frac 18}$. Since the bound for $V_N$ is asymptotically sharp (we are in the symmetric
case) there exists some $t_N$ such that $$\P_\mu\left(\frac{1}{t_N} \, \int_0^{t_N} V_N(X_s) ds \,
- \, \int V_N d\mu  \, \geq \, R\right) \, \geq \, \exp \left( - \frac 12 \, t_N \, H^*_N(R+\int
V_N d\mu)\right) \, ,$$ from which it is easy to deduce (taking $R=R'+\int V - \int V_N$ and using
that $V\geq V_N$), $$\P_\mu\left(\frac{1}{t_N} \, \int_0^{t_N} V(X_s) ds \, - \, \int V d\mu  \,
\geq \, R'\right) \, \geq \, \exp \left( - \frac 12 \, t_N \, N^{- \, \frac 18} (R'+\int V d\mu)
\right) \, .$$ Hence we have no asymptotic exponential bound.
\end{remark}
\smallskip

Since Theorem \ref{thmfsob} is not satisfactory when $F(0)>-\infty$, we shall complete it, at least
when a Poincar\'e inequality also holds.
\begin{theorem}\label{thmfsob2}
Let $F$ and $\mu$ be as in Theorem \ref{thmfsob}. Assume in addition that $\mu$ satisfies some
Poincar\'e inequality with constant $C_P$. Let $V$ such that $\int V d\mu = 0$ and $\int V^2 d\mu =
m^2 <+\infty$. Assume that $\int F^{-1}(\lambda_0 V^+) \, d\mu \, < \, + \infty$ for some
$\lambda_0 > 0$ and define $$\lambda_1 \, = \, \sup \left\{0<\lambda \leq \lambda_0 \, ; \,
F\left(\int F^{-1}(2 \lambda V \, \BBone_{\lambda V > 1/4C_P}) \, d\mu\right) \, \leq 1/4 C_P
\right\} \, .$$ Then for all $R>0$ $$\P_\nu\left(\frac{1}{t} \, \int_0^{t} V(X_s) ds \, \geq \,
R\right)\, \leq \,
\parallel \frac{d\nu}{d\mu}\parallel_{\L^2(\mu)} \, \exp \, \left\{ - t \, \sup_{0\leq \lambda \leq \lambda_1} \,
\left(R \lambda \, - \, 8 m^2 C_P \lambda^2\right)\right\} \, .$$
\end{theorem}
\begin{proof}
Using the notation in Proposition \ref{propPoinc} we  have
\begin{eqnarray}\label{eqfsob2}
\Lambda(\lambda V) & \leq & \sup_{\varepsilon \geq 0} \, \left(\frac{\varepsilon}{1+\varepsilon^2}
\, \sup \, \left\{2\lambda \, \int Vg d\mu + \varepsilon \, \left(\int \lambda V g^2  d\mu \, + \,
<Lg,g>_\mu\right) \right\}\right) \, .
\end{eqnarray}
First $\int \lambda V g d\mu \, \leq \, m \lambda $ . Next, the second term is splitting into the
sum of
\begin{eqnarray*}
\int \lambda V g^2 \, \BBone_{\lambda V \leq G_P/4} d\mu \, + \, \frac 12 \, <Lg,g>_\mu & \leq & -
(G_P/4) \, ,
\end{eqnarray*}
according to Poincar\'e, and of
\begin{eqnarray*}
\int \lambda V g^2 \, \BBone_{\lambda V > G_P/4} d\mu \, + \, \frac 12 \, <Lg,g>_\mu & \leq & \frac
12 \, F\left(\int F^{-1}(2 \lambda V \, \BBone_{\lambda V > G_P/4}) \, d\mu\right) \, ,
\end{eqnarray*}
according to the proof of Theorem \ref{thmfsob}. Using the definition of $\lambda_1$ we finally see
that for $\lambda \leq \lambda_1$, $$\Lambda(\lambda V) \, \leq \, \sup_{\varepsilon} \left(2
\varepsilon m \lambda \, - \, (1/8C_P) \varepsilon^2\right) \, \leq \, 8 m^2 \lambda^2 C_P \, .$$
The result follows.
\end{proof}

\begin{remark}
1) The existence of $\lambda_1$ is ensured by the properties of $F$ and the existence of
$\lambda_0$, while the existence of the variance of $V$ is ensured by the existence of $\lambda_0$.
Once again we obtain a Gaussian bound for small $R$ and an exponential one for large $R$.

2) As shown by Aida (\cite{Aid98}) a $F$-Sobolev inequality together with a weak Poincar\'e
inequality imply the ordinary Poincar\'e inequality. Since any absolutely continuous measure $\mu$
on a manifold with bounded from below Ricci curvature satisfies a weak Poincar\'e inequality, the
Poincar\'e inequality is automatically satisfied in this case. In particular Theorem \ref{thmfsob2}
completes the picture for the $F_\alpha$ introduced in example \ref{exF}.

3) One can easily obtain a very rough bound for $\lambda_1$. For instance if $F^{-1}=\exp
(x^\theta)$ for some $\theta>1$, the following $$\int \, \BBone_{\lambda_1 V > 1/4C_P} \, e^{(2
\lambda_1 \, V^+)^\theta} d\mu \, \leq \, e^{(1/4C_P)^\theta} \, - \, 1$$ is a sufficient
condition. Choosing $(2 \lambda_1)^\theta \leq \frac 12 \, (\lambda_0)^\theta$ and applying
Cauchy-Schwarz we obtain $$\mu(\lambda_1 V > 1/4C_P) \, \leq \, \frac{ \left(e^{(1/4C_P)^\theta} \,
- \, 1\right)^2}{\int e^{(\lambda_0 V^+)^\theta} d\mu}$$ and the left hand side is less than
$\left(\int e^{(\lambda_0 V^+)^\theta} d\mu\right) \, e^{-(\lambda_0/4C_P \lambda_1)^\theta}$
yielding an explicit condition for $\lambda_1$.
\end{remark}

\begin{remark}
In all explicit cases we considered, $x \mapsto x F(x)=G(x)$ is convex. Hence, using the
$F$-Sobolev inequality, $\Lambda(\lambda V) \leq \sup \left(\int \lambda V h d\mu - \int G(h)
d\mu\right)$ where the supremum is taken over all nonnegative $h$ such that $\int h d\mu =1$. This
kind of maximisation problem is well known in convex analysis since it relies on the calculation of
the Fenchel-Legendre transform of an integral which is a convex functional. If we relax both
constraints on $h$, one expects that this supremum is equal to $\int G^*(\lambda V) d\mu$. Actually
the situation is a little bit more intricate since $\mathbb L^1$ is not reflexive (see
\cite{Roc1,Roc2,Leo1}). Nevertheless the result we get using this (potential) bound is not
interesting. Indeed consider $G(x)=x \log(x)$ ($G(x)=+\infty$ if $x<0$) so that $G^*(u)=e^{u-1}$.
Since $G^*(0)\neq 0$ we do not obtain any interesting bound for small $R$. For instance if
$V=\BBone_A - \BBone_{A^c}$ for some $A$ with $\mu(A)=1/2$ we obtain $H(R)=R \arg \sinh (eR) -
\sqrt{1+e^2 R^2}/e$ which is negative for small $R$.

It seems that the maximisation problem taking into account the (non linear) constraints on $h$ is
not easy and we did not find any reference on it (see however \cite{Leo2} for connected problems
with linear constraints).
\end{remark}
\smallskip

Nevertheless something can be made for still stronger $F$-Sobolev inequalities. First we introduce
some definitions.

\begin{definition}\label{defyoung}
We shall say that $F$ is a contractive function if $F(x) \to +\infty$ when $x$ goes to $+\infty$
and $x \mapsto x F(x) := G(x)$ is a normalised Young function. This means that $G$ (defined on
$\R^+$) is convex, non-decreasing, satisfies $G(0)=0$ and $G(1)+G^*(1)=1$ where $G^*$ is the
Fenchel-Legendre conjugate of $G$. We shall denote by $N_G$ the corresponding gauge norm (i.e.
$N_G(f)=\inf \{u >0 \, ; \, \int G(f/u) d\mu \leq G(1)\}$).
\end{definition}

\begin{definition}\label{defFsobfort}
Let $F$ be a contractive function. \begin{itemize} \item We shall say that $\mu$ satisfies the
strong $F$-Sobolev inequality with constant $C_{SF}$ if for all $g\in D(L)$ such that $\int g
d\mu=0$ and $\int g^2 d\mu =1$ it holds $$\int g^2 F(g^2) d\mu \, \leq \, - C_{SF} \, <Lg,g>_\mu \,
.$$ \item We shall say that $\mu$ satisfies the (defective) $F$-Sobolev inequality with constants
$C_F$ and $C_b$ if for all $f\in D(L)$ such that $\int f^2 d\mu =1$ it holds $$\int f^2 F(f^2) d\mu
\, \leq \, - C_F \, <Lf,f>_\mu \, + \, C_b \, F(1) \, .$$  If $C_b=1$ we say that the inequality is
tight.
\end{itemize}
\end{definition}

Before stating the results we have in mind in the above situation, we shall discuss Definition
\ref{defFsobfort} and give some examples.

\begin{example}\label{exfsobfort}
1) It is immediate that a $F$-Sobolev inequality together with a Poincar\'e inequality imply a
strong $F$-Sobolev inequality with $C_{SF}=C_F + C_b F(1) C_P$.

2) When $F(x)=x^p$ the local version of the strong $F$-Sobolev inequality is quite useful for
studying regularising effects in p.d.e. theory for elliptic degenerate operators. In particular if
$L$ is a sub-elliptic operator in $\R^d$ satisfying some degeneracy conditions, we may associate a
natural distance to $L$ and the balls corresponding to the distance. If $d\mu/dx$ belongs to some
appropriate Muckenhoupt space, then $\mu$ will satisfy the strong $x^p$-Sobolev for some
appropriate $p$ in all balls. For precise results in this direction see e.g. Franchi \cite{Fra}
Theorem 4.5 or Lu \cite{Lu} Theorem B.

3) The defective $F$-Sobolev inequality is much more well known. First we are using ``tight''
following Bakry, while we used ``additive'' in \cite{BCR1}. Tight means that we have an equality
for $f\equiv 1$. If a non tight inequality holds (i.e. replacing $F(1)$ by a larger constant)
together with a Poincar\'e inequality, then modifying the constant $C_F$ we may obtain a tight one.
Actually this result is not proved in full generality but is proved in \cite{bakry} for $F(x)=x^p$
(the usual Sobolev inequality). Hence on a Riemannian manifold with Ricci curvature bounded from
below, we may apply Aida's result formerly recalled.

It is well known in the symmetric case (see \cite{Dav} Corollary 2.4.3) that a $x^p$-Sobolev
inequality for $p>2$ is equivalent to the ultracontractive bound $\parallel P_t f\parallel_\infty
\leq c \, t^{-s} \, \parallel f \parallel_2$ for $0<t\leq 1$ with $s=p/2(p-2)$. We shall use this
bound in the next section.

More generally $F$-Sobolev inequalities are related to super-Poincar\'e inequalities. A precise
discussion is done in \cite{w00} (also see \cite{GWa02}). In particular it is shown therein that
for $\beta \geq 1$ the measure $\mu_\beta(dx)=\frac{\exp ( - |x|^\beta)}{Z_\beta} \, dx$ satisfies
a defective $F_\alpha$-Sobolev inequality with $F_\alpha(x)=\left(\log(1+x)\right)^\alpha$ and
$\alpha=2(1-1/\beta)$. According to 1), $\mu_\beta$ thus satisfies the strong $F_\alpha$-Sobolev
inequality (since it satisfies Poincar\'e).

4) In the previous points we did not take care on the normalisation assumption for $G$. It is known
that if $G$ is not normalised one can find some $k$ such that $G(k x)$ is normalised. If $G$ is
moderate (i.e. $G(2x) \leq c G(x)$ for some $c$ and all $x$) we may replace it by its normalised
equivalent, up to a change in $C_{SF}$. Hence in the cases we discussed before, the normalisation
hypothesis is not really relevant.
\end{example}
\medskip

We conclude this section with the analogue of Theorem \ref{thmfsob2} in the contractive situation.

\begin{theorem}\label{thmfcontract}
Assume that $\mu$ satisfies a strong $F$-Sobolev inequality for some contractive function $F$. Let
$V$ such that $\int V d\mu = 0$ and $\int V^2 d\mu = m^2 <+\infty$. Then for all $R>0$
$$\P_\nu\left(\frac{1}{t} \, \int_0^{t} V(X_s) ds \, \geq \, R\right)\, \leq \,
\parallel \frac{d\nu}{d\mu}\parallel_{\L^2(\mu)} \, \exp \, \left\{ - t \, \sup_{0\leq \lambda \leq \lambda_2} \,
\left(R \lambda \, - \, (2 m^2 C_{SF}/G(1)) \lambda^2\right)\right\} \, ,$$ where
$\lambda_2=G(1)/(2 C_{SF} N_{G^*}(V))$. In particular for this bound to be interesting one needs
$N_{G^*}(V)<+\infty$.
\end{theorem}

\begin{proof}
Recall that, using the notation in Proposition \ref{propPoinc} we  have
\begin{eqnarray*}
\Lambda(\lambda V) & \leq & \sup_{\varepsilon \geq 0} \, \left(\frac{\varepsilon}{1+\varepsilon^2}
\, \sup \, \left\{2\lambda \, \int Vg d\mu + \varepsilon \, \left(\int \lambda V g^2  d\mu \, + \,
<Lg,g>_\mu\right) \right\}\right) \, .
\end{eqnarray*}
But according to the H\"{o}lder-Orlicz inequality (for normalised Young functions) $$1 \, = \, \int
g^2 d\mu \, \leq \, N_G(g^2) \, N_{G^*}(1) \, = \, N_G(g^2) \, .$$ It follows that $\int  G(g^2)
d\mu \geq G(1) \, N_G(g^2)$. Using H\"{o}lder-Orlicz inequality again and the strong $F$-Sobolev
inequality, it holds
\begin{eqnarray*}
\int \lambda V g^2  d\mu \, + \, <Lg,g>_\mu & \leq & N_{G^*}(\lambda V) \, N_G(g^2) \, - \,
(1/C_{SF}) \, \int G(g^2) d\mu \\ & \leq & \left(\lambda \, N_{G^*}(V) \, - \, (G(1)/C_{SF})\right)
\, N_G(g^2) \\ & \leq & \lambda \, N_{G^*}(V) \, - \, (G(1)/C_{SF}) \, ,
\end{eqnarray*}
provided $\left((G(1)/C_{SF}) - \lambda \, N_{G^*}(V)\right) >0$. It follows that $$\Lambda(\lambda
V) \, \leq \, \sup_{\varepsilon \geq 0} \left\{2 m \lambda \varepsilon \, - \, \left((G(1)/C_{SF})
- \lambda \, N_{G^*}(V)\right) \, \varepsilon^2\right\} \leq m^2 \lambda^2 /\left((G(1)/C_{SF}) -
\lambda \, N_{G^*}(V)\right) \, .$$  The proof is completed.
\end{proof}

The previous proof is certainly simpler than the one of Theorem \ref{thmfsob2}, while both theorems
apply to similar measures, for instance the $\mu_\beta$'s for $\beta > 1$. It is quite difficult to
compare the bounds in both Theorems on this example, but the one in Theorem \ref{thmfcontract} has
to be worse in general since it lies on the rough use of H\"{o}lder and Orlicz norms.
\smallskip

\begin{remark}
We have not discussed here the use of another type of functional inequalities called transportation cost inequalities (in path space), namely
$$\forall\nu,\qquad W_p(\nu,\mu)\le\sqrt{2C~\int\log(d\nu/d\mu)d\nu}$$
where $W_p$ is the usual Wasserstein distance, leading to Gaussian type of deviation inequalities for Lipschitz test function $V$. But the proofs are very different in spirit as they rely on the verification of some square exponential integrability for some norm on the path space and a (dependent) tensorization property. We refer to Djellout-Guillin-Wu \cite{DGW}. Note that the results obtained there are reminiscent of an assumption of a logarithmic Sobolev inequality to hold. However they do not rely on the knowledge of the invariant measure but on the conditions on the drift and diffusion coefficient. In the same spirit, one may also use Poincar\'e inequalities or logarithmic Sobolev inequalities on path space (see \cite{CHL} for example) combined with Herbst's argument, but they are much more difficult to prove and does not give good bound for large time asymptotic.
\end{remark} 

\bigskip

\section{\bf Polynomial and sub-exponential bounds}\label{secpoly}

\subsection{\bf The case of bounded
$V$'s.}\label{subsecbounded}

In this subsection we shall assume that the semi-group $P_t$ satisfies the following decay
property:
\begin{assumption}\label{dec}
there exists some non increasing function $\eta$ defined on $[0,+\infty[$ such that for all bounded
$f$ and all $t$, $$\Var_\mu(P_t f) \, \leq \, \eta(t) \, \parallel f - \int f d\mu
\parallel_\infty^2 \, .$$
\end{assumption}

It is known (see \cite{r-w} Theorem 2.1 and 2.3) that Assumption \ref{dec} is formally equivalent
to a weak Poincar\'e inequality (WPI). More precisely, if $\mu$ satisfies a (WPI), i.e. for all
$s>0$ and all bounded $g$,
\begin{equation}\label{wpi}
\Var_\mu(g) \, \leq \, - \, \beta(s) \, <Lg,g>_\mu \, + \, s \, \parallel g - \int g d\mu
\parallel_\infty^2 \, ,
\end{equation}
for some non increasing $\beta$, then \eqref{dec} holds with $$\eta(t) \, \leq \, 2 \, \inf
\left\{s>0 \, ; \, \beta(s) \, \log(1/s) \, \leq \, 2t \right\} \, .$$ Conversely, in the symmetric
case (ore more generally if $L$ is a normal operator), if $\eta$ is decreasing with inverse
function $\eta^{-1}$, then \ref{dec} implies a (WPI) with $$\beta(t) \, = \, 2t \, \inf_{s>0}
\left( \frac 1s \, \eta^{-1}(s \, \exp[1-s/t])\right) \, .$$ In particular, if $\eta(t) \leq
e^{-\delta t}$ for some $\delta >0$, Assumption \ref{dec} implies the (true) Poincar\'e inequality. In the
sequel we may thus assume that $\eta$ is decaying slower than an exponential. We shall give
explicit examples later.
\smallskip

Assumption \ref{dec} is clearly connected to \underline{mixing} properties. Indeed recall the

\begin{definition}\label{defmix}
The strong mixing coefficient $\alpha(r)$ is defined as $\alpha(r) =  \sup_{s,F,G} \{|Cov(F,G)| \}$
where $F$ (resp. $G$) is $\mathcal F_s$ (resp. $\mathcal F_{s+r}$) measurable, non-negative and
bounded by 1.
\end{definition}

We then have

\begin{proposition}\label{propmix}
If Assumption \ref{dec} holds then the stationary process is strongly mixing with $\alpha(r)\leq
\sqrt{\eta(r)}$. If in addition $\mu$ is symmetric we may choose $\alpha(r)\leq \eta(r/2)$.

Conversely, if the stationary process is strongly mixing, then Assumption \ref{dec} holds with
$\eta(r) \leq \alpha(r)$.
\end{proposition}
\begin{proof}
If $F$ and $G$ are centred and bounded by 1, we may apply the Markov property to get $$\E_\mu[FG]
= \E_\mu [F \, \E[G/X_{s+r}]] = \E_\mu[F \, P_r g(X_s)] $$ where $g$ is centred and bounded by 1.
Hence $$|\E_\mu[FG]| \leq \E_\mu[|P_r g(X_s)|] = \int |P_r g| d\mu \leq \sqrt{\eta(r)} \, .$$ In
the symmetric case $$\E_\mu[F \, P_r g(X_s)] = \E_\mu[F(X_{s-.}) \, P_r g(X_0)]=\E_\mu[f(X_0) \,
P_r g(X_0)]= \int P_{r/2}f \, P_{r/2}g \, d\mu$$ and we conclude using Cauchy Schwarz again.

For the converse, taking $F=P_r f(X_0)$ and $G=f(X_r)$ for $f$ centred and bounded by one
furnishes the result.
\end{proof}

The point is that moment bounds for sums of strongly mixing sequences (extending Rosenthal's
inequalities in the independent case) are known. A large part of them are due to Doukhan and his
coauthors and may be found in Doukhan's book \cite{douk}. However we found the most refined version
we shall use in Rio's book \cite{rio}.

\begin{proposition}\label{proppolyn}
Assume that $\mu$ satisfies Assumption \ref{dec} for some $\eta$ satisfying for some integer $k$,
$M_k(\alpha) := \sup_r (1+r)^k \, \alpha(r) < +\infty$ with $\alpha$ as in Proposition
\ref{propmix}.

Then there exists a constant $C(k)$ such that for all $V$ with $\sup |V|= 1$, all $0<R\leq 1$ and
all $t<[t]/(1-R)$ where $[t]$ is the integer part of $t$
\begin{equation}\label{eqpolyn}
\P_\mu\left(\frac 1t \, \int_0^t V(X_s) ds \, - \, \int V d\mu  \, \geq \, R\right) \, \leq \,
\frac{ C(k) M_k(\alpha)}{ t^k \, \left(R-(1-([t]/t))\right)^{2k}} \, .
\end{equation}
\end{proposition}
\begin{proof}
Denote by $Y_j=\int_{j-1}^j \, V(X_s) ds - \int V d\mu$. Then $Y_j$ is a ($\P_\mu$) stationary
sequence of strongly mixing centred random variables with mixing coefficient $\alpha(r-1)$. Thanks
to our hypothesis on $\alpha$ we may apply Theorem 2.2 in \cite{rio} (see (2.23) p.40 therein)
which yields $\E_\mu [(\sum_1^n \, Y_j)^{2k}] \, \leq \, C(k) \, M_k(\alpha) \, n^k \, .$ The
result follows by using Markov inequality and the fact that $V$ is bounded by 1.
\end{proof}

In the previous result one can obtain explicit bounds for $C(k)M_k(\alpha)$ as shown by Doukhan and
Portal (see \cite{douk} chapter 1.4).

In the examples below again, $-<Lf,f>_\mu = 1/2 \, \int |\nabla f|^2 d\mu$.
\begin{example}
1) If $\mu(dx)=c (1+|x|)^{-(d+p)} dx$ ($p>0$) on $\R^d$ it is shown in \cite{BCR2} that (WPI) holds
with $\beta(s)=c(d) s^{-2/p}$. Actually this result is shown for $d=1$ but extends to the $d$
dimensional case since the tensorized $1$-dimensional measure is equivalent to the $d$-dimensional
one (of course all constants depend on $d$). Hence we may choose $\eta(t)=c(d,p) (\log t /
t)^{p/2}$. The bound in Proposition \ref{proppolyn} is thus available for $p>2k$ in the symmetric
case, and $p>4k$ in the non symmetric-one.

2) If $\mu(dx)=c e^{-|x|^p} dx$ for some $1\geq p>0$, we obtain similarly $\eta(t)=c(d,p) \, e^{-
c' t^{\frac{p}{2-p}}}$. We can thus obtain any polynomial bound. Of course, in this case one
expects a better bound. We shall see how to get such a bound below.
\end{example}

In order to get sub-exponential bounds, we recall the following moment inequality from \cite{rio}
Theorem 2.5, that holds for a ($\P_\mu$) stationary sequence of strongly mixing centred random
variables $Y_j$ bounded by 1
\begin{equation}\label{riofin}
\E_\mu [|\sum_1^n \, Y_j|^{2k}] \, \leq \, (4nk)^k \, \int_0^1 \, (\alpha^{-1}(u)\wedge n)^k \, du
 \, .
\end{equation}
Note that \eqref{riofin} allows us to give an explicit bound for $C(k)$ in Proposition \ref{proppolyn}, but
with a slightly worse speed. Indeed if $\alpha(n)=c n^{-k}$ we get $$\E_\mu [|\sum_1^n \,
Y_j|^{2k}] \, \leq \, (4nk)^k \, c \, (1+ \log c + k \log n)$$ recovering \eqref{eqpolyn} with an
extra logarithm.

Recall now the following elementary $$\limsup_{q\to +\infty} q^{-1} \left(\int_0^1 \, \log^q(1/u)
\, du\right)^{1/q} \, \leq \, 1/e \, .$$ If $\alpha(n)=c \, e^{- c' n^{p/(2-p)}}$ it follows that
there exists some $k_0$ depending only on $c$ such that $$\int_0^1 \, (\alpha^{-1}(u)\wedge n)^k \,
du \, \leq \, \left(\frac{(2 - p)k}{e \, p \, c'}\right)^{\frac{(2 - p)k}{p}} \, ,$$ for all $k\geq
k_0$. Hence $\E_\mu [|\sum_1^n \, Y_j|^{2k}] \, \leq \, \left( 4n \left(\frac{2-p}{e \, p \,
c'}\right)^{\frac{2-p}{p}} \, k^{2/p}\right)^k $.

Using Markov inequality it thus holds
\begin{equation*}
\P_\mu\left(\int_0^n V(X_s) ds \, - \, \int V d\mu  \, \geq \, S \sqrt{n}\right) \, \leq \, e^{-
\frac{2-p}{p} \, k} \, \left(2 \left(\frac{2-p}{p \, c'}\right)^{\frac{2-p}{2p}} \, k^{1/p} \,
(1/S) \right)^{2k} \, .
\end{equation*}
We then choose $k^{1/p}=(pc'/2-p)^{2-p/2p} \, S/2$ provided it is greater than $k_0^p$. Finally we
have obtained,
\begin{proposition}\label{propsubexp}
Assume that $\mu$ satisfies Assumption \ref{dec} with $\alpha(s)=c \, e^{- c' s^{p/(2-p)}}$ as in
Proposition \ref{propmix}, for some $0<p\leq 1$.

Then there exists a constant $k_0$ depending on $c$, such that for all $V$ with $\sup |V|= 1$, all
$0<R\leq 1$ and all $t<[t]/(1-R)$ where $[t]$ is the integer part of $t$
\begin{equation}\label{eqsubexp}
\P_\mu\left(\frac 1t \, \int_0^t V(X_s) ds \, - \, \int V d\mu  \, \geq \, R\right) \, \leq \, \exp
\left\{
 - \, c(p) \, \left((R-(1-([t]/t))) \sqrt t\right)^p \right\} \, ,
\end{equation}
with $c(p) = \frac{2-p}{p} \, (1/2)^p \, (pc'/2-p)^{2-p/2}$, provided $(R-(1-([t]/t))) \sqrt t \,
\geq \, 2 k_0^{1/p} \, (2-p/pc')^{2-p/2p}$.
\end{proposition}
\medskip

\begin{remark}
If the previous result is in accordance with the C.L.T. (that holds as soon as $\int_0^\infty
\alpha(s) ds < +\infty$), it is of course worse than the ones we obtained in the first section.
Indeed for $p=1$ we recover a convergence rate $e^{-C \sqrt t}$ (for some fixed $R$) while we know
that (at least in the symmetric case) a Poincar\'e inequality holds, hence Proposition
\ref{propPoinc} gives a convergence rate $e^{-C t}$.

This fact suggests that we may loose something in the time discretization. At the same time we may
ask whether it is possible to use the semi-group structure to calculate $$G_{2k}(t) :=
\E_\mu[|\int_0^t V(X_s) ds|^{2k}]$$ or not. If $k\in \N$ it is possible to study the variations of
$G_{2k}$, at least in the symmetric diffusion case. We assume now that $\int V d\mu = 0$. Then
$$G'_2(t)= 2 \E_\mu[V(X_t) \, \int_0^t V(X_s) ds]= 2 \int_0^t \int (P_{s/2}V)^2 d\mu ds \leq 2
\int_0^t \eta(s) ds \, .$$ Hence if $\eta \in \L^1(\R^+,dt)$ we obtain that $G_2(t) \leq c_2 t$.

Using integration by parts and symmetry one can show that
\begin{eqnarray*}
(G_4)''(t)&=& 12 \, \E_\mu\left[ V(X_0) \, V(X_t) \, \left(\int_0^t V(X_s) ds\right)^2\right] \,\\
&=& 24 \, \E_\mu\left[ V(X_0) \, V(X_t) \, \left(\int_{0}^{t/2} V(X_s) ds\right)^2\right] \, + \\ &
& + \, 24 \, \E_\mu\left[ V(X_0) \, V(X_t) \, \left(\int_{t/2}^{t} V(X_s)
ds\right)\left(\int_{0}^{t/2} V(X_s) ds\right)\right]
\end{eqnarray*}
The first term in the above sum can be bounded by $6 t^2 \, \sqrt{\eta(t/2)}$, so that if this last
quantity is in $\L^1(dt)$ it furnishes a contribution $c_4 t$ again to $G_4$. The second term in
the sum can be rewritten with the help of the function $$H(x)=\E_x \left[\int_0^{t/2} V(X_{t/2})
V(X_s) ds\right]=\int_0^{t/2} P_s\left(V P_{t/2 -s}V\right)(x) ds \, .$$ It yields
$$\E_\mu\left[ V(X_0) \, V(X_t) \, \left(\int_{t/2}^{t} V(X_s)
ds\right)\left(\int_{0}^{t/2} V(X_s) ds\right)\right] = $$
\begin{eqnarray*}
& = & \E_\mu\left[V(X_0) \left(\int_{0}^{t/2} V(X_s) ds\right) H(X_{t/2})\right] \\ & = &
\E_\mu\left[H(X_0) V(X_{t/2}) \left(\int_{t/2}^{t} V(X_s) ds\right)\right] \, \left(= \int H^2 d\mu\right) \, \\
& = & \int_0^{t/2} \int V \, P_s(V P_{t/2 -s}H) \, d\mu \, ds =  \int_0^{t/2} \int (P_sV) \, V \,
(P_{t/2 -s}H) \, d\mu \, ds \\ & = & \int_0^{t/2} \int (P_sV) \, V \, (P_{t/2 -s}(H-\int H d\mu))
\, d\mu \, ds \, + \, (\int H d\mu) \, \int_0^{t/2} \int (P_sV) \, V \, d\mu \, ds
\\& \leq & \int_0^{t/2} \sqrt{\eta(s)} \, \sqrt{\eta(t/2 - s)} \, \Var_\mu^{1/2}(H) \, ds \, +
 \, (\int H d\mu) \, \int_0^{t/2} \eta(s/2) ds \, .
\end{eqnarray*}
But $$\int H d\mu = \int_0^{t/2} \int \, V \, P_{t/2 -s}V  \, d\mu \, ds \, \leq \, \int_0^{t/2}
\eta(s/2) ds$$ is assumed to be bounded by $d_4$ (see the control of $G_2$). Since $\Var_\mu(H)
\leq \int H^2 d\mu$ it follows that $$\Var_\mu^{1/2}(H) \leq \int_0^{t/2} \sqrt{\eta(s)} \,
\sqrt{\eta(t/2 - s)} ds + d_4 \leq (t/2) \sqrt{\eta(t/4)} + d_4 \, .$$ Since we formerly assumed
that $t^2 \sqrt{\eta(t/2)}$ goes to 0, the Variance is bounded and consequently so is $(G_4)''$
yielding a bound $c_4 t^2$ for $G_4$.

Unfortunately, it seems difficult to iterate the procedure and to get explicit expressions for the
constants. Furthermore, one suspects that a clever study will yield $G_{2k}(t) \leq c_{2k} t^k$,
that is the same behaviour as in the discrete case. It does not seem necessary to go further.
\end{remark}

\medskip

\subsection{Unbounded $V$'s.}

If $V$ is no more bounded, or does not fulfill the hypotheses of one of the result in the second
section, one can get some bound by truncating $V$. We shall briefly indicate how to do on an
example.

For instance for a centred $V$ such that $\int |V| d\mu < +\infty$, and all $K>0$
$$\P_\mu\left(\frac 1t \, \int_0^t V(X_s) ds \, \geq \, R\right) \leq
$$
\begin{eqnarray*}
& \leq & \P_\mu\left(\frac 1t \, \int_0^t (V\wedge K\vee -K)(X_s) ds \, - \, \int (V\wedge K\vee
-K) d\mu \, \geq \, R/2 - \int (V\wedge K\vee -K) d\mu\right) + \\ & & + \P_\mu\left(\frac 1t
\int_0^t |V|\BBone_{|V|\geq K}(X_s) ds \, \geq \, R/2\right) \, = \, A + B.
\end{eqnarray*}

If $\int |V|^S d\mu < +\infty$, $B$ can be bounded by $K^{m-S} \, R^{-m} \, 2^m \, \E_\mu[|V|^S]$
for all $1\leq m<S$. If $\int e^{u|V|} d\mu < +\infty$ for some $u>0$, we have $B \leq e^{-\lambda}
\, \int e^{\lambda |V| \sqrt{2/KR}} d\mu$ as soon as $\lambda \sqrt{2/KR} \leq u$ (just summing up
the previous bounds for $m=S/2$).
\smallskip

In order to obtain a bound for $A$ we may use the appropriate results in section \ref{Intro} or in
the previous subsection.

If we assume for example that $\mu$ satisfies a Poincar\'e inequality, and that $K$ is such that
$R/4 \geq \E_\mu[|V|^S]/K^{S-m}$ we may apply Proposition \ref{propPoinc} and obtain a bound for
$A$ in the form $exp \{- t R^2 / 128 C_P \, K^2\}$.

Choosing $m=S/2$ (provided $S\geq 2$), it is not difficult to see that the (almost) optimal choice
is given by $K = c R \sqrt{t} / \sqrt{\log(2 tR^4/S)}$ with $t$ large enough for this expression to
be meaningful and the previous constraint between $R$ and $K$ to be satisfied.

We thus obtain a bound $$ C(S) \, \log^{S/4}(2 tR^4/S) \, R^{-S/2} \, t^{-S/4} \, ,$$ for $t$ large
enough, $\mu$ satisfying Poincar\'e and $\E_\mu[|V|^S]<+\infty$ for some $S\geq 2$.

In the same way, if $\int e^{u|V|} d\mu < +\infty$, we first choose $\lambda=u \sqrt{KR/2}$, and a
similar method yields a bound $$C(u) e^{-c(u) t^{1/5} R^{4/5}} \, ,$$ for $t$ large enough.
\bigskip

\section{\bf About the initial measure}\label{initial}

In this final section we shall see what can be said for the initial measure $\nu$. As for the
latter subsection, we shall not state general results, but give some hints in various situations.
Of course we shall discuss how to get deviation bounds for $\P_\nu(\mathcal F)$ which are not
simply given by $\P_\mu^{1/p}(\mathcal F) \,
\parallel d\nu/d\mu\parallel_q$.
\medskip

{\bf A.} \quad We have seen in sections \ref{sob} and \ref{Intro} that we may take some initial
measure $\nu$ such that $d\nu/d\mu \in \L^2(\mu)$. As remarked by Wu \cite{Wu0} p.441-442, we may
replace this assumption by $d\nu/d\mu \in \L^q(\mu)$ for $1\leq q < +\infty$, provided we replace
$\Lambda$ by
$$\Lambda_p(V) := \sup \, \left\{ \int \, V|f|^p \, d\mu \, + \, <\textrm{sgn}(f)|f|^{p-1},Lf>_\mu \, ; \, f\in D_p(L) \,
\textrm{ and } \, \int |f|^p d\mu =1 \right\} \, ,$$ where $p$ and $q$ are conjugate. If $L$ admits
a carr\'e du champ $\Gamma$, one can integrate by parts and get $$<\textrm{sgn}(f)|f|^{p-1},Lf>_\mu
\, = \, - \, (4(p-1)/p^2) \, \int \Gamma(|f|^{p/2}) \, d\mu \, ,$$ so that defining $g=|f|^{p/2}$
we obtain that
$$\Lambda_p(V) \, = \, (4(p-1)/p^2) \, \Lambda( (p^2/4(p-1)) V)$$ at least for a bounded $V$ (remark
that $(p-1)/p^2=(q-1)/q^2$).

Hence all the results in sections \ref{sob} and \ref{Intro} are still true, up to the constants,
for $1<q<+\infty$. For instance we get an additional  constant $4(p-1)/p^2$ in Proposition
\ref{propPoinc}. Since the interesting $q$'s are less than 2, the interesting $p$'s are greater
than 2 and this bound is better than the $1/p$ obtained via H\"{o}lder.
\medskip

{\bf B.} \quad If $\mu$ is symmetric we may argue as follows : let $\mathcal A$ be a $\sigma (X_s,
u\leq s \leq t)$ measurable subset and denote by $R_t$ the time reversal at time $t$. Then
\begin{eqnarray*}
\P_\nu(\mathcal A) & = & \E_\mu \left[\frac{d\nu}{d\mu}(X_0) \, \BBone_{\mathcal A}\right] = \E_\mu
\left[\frac{d\nu}{d\mu}(X_t) \, \BBone_{\mathcal A}\circ R_t\right] = \E_\mu \left[\left(P_u
\frac{d\nu}{d\mu}\right)(X_{t-u}) \, \BBone_{\mathcal A}\circ R_t\right] .
\end{eqnarray*}
If $V$ is centred and bounded by 1, the set $\left\{\frac 1t \, \int_0^t V(X_s) ds \geq R\right\}$
is included in $$\mathcal A := \left\{\frac 1t \, \int_u^t V(X_s) ds \geq (R-(u/t))\right\}$$ to
which we may apply the previous trick.

In particular, if the semi-group is ultracontractive (i.e. there exists some $u>0$ such that $P_u$
is mapping continuously $\L^1(\mu)$ in $\L^\infty(\mu)$) we obtain a nice bound. Notice that $P_v$
is also mapping continuously $\L^1$ in $\L^2$ for some $v\leq u$, so that using reversibility again
we may directly use section \ref{Intro}, with a possible better constant.
\smallskip

If the semi-group is only hypercontractive, i.e. if $\mu$ satisfies some log-Sobolev inequality, we
know that relative entropy is exponentially decaying. Denote by $H(h):= \int h \log h d\mu$ for any
density of probability, and by $h_\nu := d\nu/d\mu$. If $H(h_\nu)<+\infty$ it holds $H(P_u h_\nu)
\leq e^{- u/C_{LS}} H(h_\nu)$. It is easily seen that $$\int \exp \left(\BBone_{\mathcal A} - (e-1)
\P_\mu(\mathcal A)\right) \, d\P_\mu \, \leq \, 1$$ so that using the variational definition of $H$
and reversibility again we get
$$\P_\nu(\mathcal A) = \E_\mu \left[P_u h_\nu(X_{u}) \, (\BBone_{\mathcal A} - (e-1)\P_\mu(\mathcal A)+(e-1)\P_\mu(\mathcal A))\right]
 \leq  H(P_u h_\nu) + (e-1) \P_\mu(\mathcal A)  \, .$$ Choosing $u=Rt/2$ we thus obtain
$$\P_\nu(\mathcal A) \leq (e-1)\, e^{- \frac{t R^2}{32 C_P \Var_\mu(V)}} + H(h_\nu) \, e^{-
\frac{tR}{2C_{LS}}} \, .$$ If the semi-group is only Orlicz-hypercontractive in the sense of
\cite{BCR1} we do not know whether it is possible to extend the argument to a little bit more
integrable initial densities or not. Indeed we did not find the ad-hoc quantity replacing relative
entropy.
\medskip

{\bf C.} Finally we shall see on a family of examples what can happen when $\nu$ is no more
absolutely continuous with respect to $\mu$. Actually we shall consider on $\R^d$ a diffusion
process $$X_t^x= x + B_t - \int_0^t \nabla U (X_s^x) ds \, ,$$ where $x\in \R^d$ and $B_.$ is a
standard Brownian motion. We shall assume that $U$ is $C^3$, and that there exists  some function
$\psi$ going to $+\infty$ when $|x| \to \infty$ so that $\frac 12 \Delta \psi - \nabla U . \nabla
\psi$ is bounded from above. These assumptions ensure the existence of an unique non explosive
strong solution. Furthermore the underlying Markov process $X_.$ is $\mu$ symmetric for $d\mu =
Z^{-1} e^{-2U} dx$ where $Z$ is a normalising constant.

For such a process it is known that the law of $X_t^x$ is absolutely continuous w.r.t. $\mu$. We
shall denote by $h^x_t$ its density.

If $|\nabla U|^2(y) - \Delta U (y) \, \geq \, - C_m > - \infty$ for all $y$, one can show that
$\int h^x_t \log_+^p (h_t^x) d\mu < +\infty$ for all $p\geq 1$ (see \cite{CGG} Proposition 5.1), so
that in particular, if the semi-group is hypercontractive (or ultracontractive) we may apply the
ideas in {\bf B}.

Actually  one can expect a much better integrability and it is shown in \cite{CGG} section 5.2 that
for $U(y)=|y|^q$ with $1\leq q \leq 2$, $h_t^x \in \L^{\infty}(\mu)$ for all $t>0$ ($U$ is not
$C^3$ but all the previous discussion is still available).

Indeed we discovered with the help of P.A. Zitt that actually, with our previous assumptions,
$h_t^x \in \L^2(\mu)$.

To prove it, as in \cite{CGG} we follow the idea of \cite{Ro99} Thm 3.2.7. Replacing the convex
$\gamma$ therein by $\gamma(y)=y^2$ we obtain $$\int (h_t^x)^2 d\mu \leq Z \, e^{2U(x)} \,
\E\left[e^{-2v(B_t)} \, e^{- \frac 12 \, \int_0^t [|\nabla U|^2 - \Delta U](B_s) ds} \right] \leq Z
\, e^{2U(x)} \, e^{\frac 12 \, C_m t} $$ where $e^{-2v(y)} = (2\pi t)^{-d/2} \, e^{- |y-x|^2/2t}$.
Hence we may directly apply the results in the first two sections.
\bigskip

\bibliographystyle{plain}
\bibliography{adfunc}

\end{document}